\documentclass[11pt]{article}

\usepackage{latexsym}
\usepackage{amsfonts}
\usepackage{amsmath}

\usepackage[dvips]{color}

\newtheorem{thrm}{Theorem}[section]
\newtheorem{lemma}[thrm]{Lemma}
\newtheorem{prop}[thrm]{Proposition}

\newtheorem{cor}[thrm]{Corollary}

\numberwithin{equation}{section}

\usepackage{enumerate}

\usepackage{amssymb}
\usepackage{mathtools}
\mathtoolsset{showonlyrefs}
\usepackage{mathrsfs}
\usepackage{comment}
\usepackage{hyperref}
\usepackage{xcolor}

\def\E{\mathbb{E} }
\def\P{\mathbb{P} }
\def\Q{\mathbb{Q} }
\def\R{\mathbb{R} }
\def\N{\mathbb{N} }

\makeatletter
\begin{document}
\allowdisplaybreaks

\title{\Large \bf{Tail probability of maximal displacement in critical branching L\'{e}vy process with stable branching}
	\footnote{The research of this project is supported
     by the National Key R\&D Program of China (No. 2020YFA0712900).}}
\author{ \bf  Haojie Hou \hspace{1mm}\hspace{1mm}
	Yiyang Jiang \hspace{1mm}\hspace{1mm}
Yan-Xia Ren\footnote{The research of this author is supported by NSFC (Grant Nos. 12071011 and 12231002) and The Fundamental
Research Funds for the Central Universities, Peking University LMEQF.
 } \hspace{1mm}\hspace{1mm} and \hspace{1mm}\hspace{1mm}
Renming Song\thanks{Research supported in part by a grant from the Simons
Foundation
(\#960480, Renming Song).}
\hspace{1mm} }
\date{}
\maketitle

\begin{abstract}

Consider a critical branching L\'{e}vy process $\{X_t, t\ge 0\}$ with branching rate $\beta>0, $ offspring distribution $\{p_k:k\geq 0\}$ and spatial motion $\{\xi_t, \Pi_x\}$.  For any $t\ge 0$, let $N_t$ be the collection of particles
alive at time $t$, and, for any $u\in N_t$, let $X_u(t)$ be the position of  $u$ at time $t$.
We study the tail probability of the maximal displacement $M:=\sup_{t>0}\sup_{u\in N_t} X_u(t)$ under the assumption
 $\lim_{n\to\infty} n^\alpha \sum_{k=n}^\infty p_k =\kappa\in(0,\infty)$ for some $\alpha\in (1,2)$,  $\Pi_0(\xi_1)=0$ and  $\Pi_0 (|\xi_1|^r)\in (0,\infty)$ for some $r> 2\alpha/(\alpha-1)$.
Our main result is a generalization of the main result of Sawyer and Fleischman (1979) for branching Brownian motions and that of Lalley and Shao (2015) for branching random walks, both of which are proved under the assumption $\sum_{k=0}^\infty k^3 p_k<\infty$.

\end{abstract}

\medskip

\noindent\textbf{AMS 2020 Mathematics Subject Classification:}  60J80;  60G40; 60G51

\medskip

\noindent\textbf{Keywords and Phrases}:
Branching  L\'{e}vy process,
critical branching process, Feynman-Kac representation.

\section{Introduction and notation}

\subsection{Introduction}

A branching random walk is a discrete-time Markov process defined as follows: at time $n=0$, there is a particle at $0\in \R$. At time $n=1$, this particle dies and splits into a finite number of offspring. The particle configurations of the offspring relative to their parent is a copy of a point process $\mathcal{L}$. At time $n=2$, the individuals alive at time $1$ repeat their parent's behavior and the process goes on.  Denote the law by $\mathbb{P}$.  We will use $N_n$ to denote the set of particle alive at time $n$ and for $u\in N_n$, the position of $u$ is denoted by $X_u(n)$.

Now we consider the special case  $\mathcal{L}= \sum_{i=1}^B \delta_{X_i}$, where $B$ is a non-negative integer valued random variable  with $\P(B=k)=p_k$ and $X_1, X_2,...$ are iid $\mathbb{Z}$-valued random variables independent of $B$ with common distribution $\{\mu_k, k\in \mathbb{Z}\}$.
We say that this process is critical if
\begin{align}
	\E(B)= \sum_{k=0}^\infty kp_k=1.
\end{align}
Since the total mass of the branching random walk is a Galton-Waston process, a critical branching random walk must extinct in finite time, which implies that the following maximal displacement $M$ is a finite random variable:
\begin{align}\label{Destrete-Maximal-def}
	M:= \sup_{n\in \N} \sup_{u\in N_n} X_u(n)
\end{align}
with the convention $\sup_{u\in N_n} X_u(n) =-\infty$ if $N_n =\emptyset.$ \cite{LS15} proved that if
\begin{align}\label{LS-assum}
	\sum_{k=0}^\infty k^3 p_k <\infty,\quad  \sum_{k\in \mathbb{Z}} k\mu_k=0,  \quad \sum_{k\in \mathbb{Z}} \left|k\right|^{4+\varepsilon} \mu_ k<\infty
\end{align}
for some $\varepsilon>0$, then
\begin{align}
	\lim_{x\to +\infty} x^2 \P\left(M\geq x\right)  =\frac{6 \eta^2}{\sigma^2 },
\end{align}
where $\eta^2:= \sum_{k\in \mathbb{Z}} k^2 \mu_k$ and $\sigma^2:= \sum_{k=0}^\infty k^2 p_k -1$.

Now we turn to the continuous time and space case, the branching L\'{e}vy process  in the sense of \cite{Kyp99}.
Let $\left(\xi_t, \Pi_x\right)$ be a L\'{e}vy process with $\xi_0=x$.
A  branching
L\'{e}vy process is defined as follows: initially there is a particle at
$x\in \R$ and it moves according to $(\xi_t, \Pi_x)$.
After an exponential time with parameter $\beta>0$,
independent of the motion, it dies and produces $k$ offspring with probability $p_k$,
$k\geq 0$.
The offspring move independently according to $\xi$ from the place where they are born and obey the same branching mechanism as their parent. Denote the law by $\P_x$ and $\P:=\P_0$.
 In this paper we focus on  the critical case, i.e., we always assume  that $\{p_k: k\geq 0\}$ satisfies $\sum_{k=0}^\infty kp_k =1$.

Similarly, we define the maximal position by
\begin{align}
	M:= \sup_{t\geq 0} \sup_{u\in N_t} X_u(t),
\end{align}
where $N_t$ is the set of particles alive at time $t$ and $X_u(t)$ is the position of $u\in N_t$.
When the spatial motion $\xi$ is a standard Brownian motion,  \cite{SF79}
proved that under the assumption $\sum_{k=0}^\infty k^3p_k<\infty$,
\begin{align}\label{Tail-probability}
	\lim_{x\to +\infty} x^2 \P\left(M\geq x\right)  =\frac{6}{\sigma^2 }
\end{align}
with $\sigma^2= \sum_{k=0}^\infty k^2 p_k -1$. \cite{Profeta22} extended \eqref{Tail-probability} to the case when  $\xi_t$ is a spectrally negative branching L\'{e}vy process and $\sum_{k=0}^\infty k^3 p_k<\infty$.
When the spatial motion is a $\gamma$-stable process
with index $\gamma\in (0,2)$, $\sum_{k=0}^\infty k^3 p_k<\infty$ and $\beta=1$,  \cite{LS16}  and  \cite{Profeta21} proved that
$$
\lim_{x\to +\infty} x^{\gamma/2} \P\left(M\geq x\right)  =\sqrt{\frac{2}{\gamma}}.
$$
For results where the spatial motion is a general spectrally negative L\'{e}vy process, see \cite{Profeta22}.

\subsection{Main result}

The main aim of this paper is to study the tail probability of $M$ when the offspring distribution $\{p_k:k\geq 0\}$ is in the domain of attraction of an $\alpha$-stable distribution with index $\alpha\in (1,2)$ and the spatial motion has light tails. Suppose that there exist constants $\kappa>0$ and $\alpha\in (1,2)$ such that
\begin{align}\label{condition-1}
	\lim_{n\to\infty} n^\alpha \sum_{k=n}^\infty p_k = \kappa.
\end{align}
Assume that
\begin{align}\label{condition-2}
	\Pi_0 (\xi_1)=0,\quad \eta^2:= \Pi_0 (\xi_1^2)\in (0, \infty).
\end{align}
Our main result is as follow:

\begin{thrm}\label{thm1}
	 If
	 \begin{align}\label{condition-3}
	 \Pi_0\left( \left| \xi_1\right| ^r\right)<\infty\quad \mbox{for some } r> \frac{2\alpha}{\alpha-1},
	 \end{align}
	then
	 \begin{align}\label{main}
	 	\lim_{x\to\infty} x ^{\frac{2}{\alpha-1}} \P \left(M\geq x\right) = \left(\frac{(\alpha+1)\eta^2 }{\beta \kappa(\alpha-1)\Gamma(2-\alpha)}\right)^{\frac{1}{\alpha-1}},
	 \end{align}
	 where $\Gamma(z):= \int_0^\infty t^{z-1}e^{-t}\mathrm{d}t$ is
the Gamma function.
\end{thrm}

Note that $\frac{2\alpha}{\alpha-1}>4$, so the spatial motion has at least finite 4-th moment.

Our argument of proving the above main result is an adaptation of that of Lalley and Shao \cite{LS15}. Our assumption \eqref{condition-1} on branching mechanism is weaker than the assumption \eqref{LS-assum} in \cite{LS15}.
Under our assumption that the spatial motion has light tails, the weaker assumption above on the branching mechanism does
not cause too much trouble.
The assumption \eqref{condition-1} only changes the behavior of $f$, defined in \eqref{Function-F} below,
from $f(v)= Cv(1+o(1))$ to $f(v)= Cv^{\alpha -1}(1+o(1))$ for some constant $C>0$. In \cite{LS15}, the  explicit solution  of the following problem
$$
\left\{\begin{array}{rl}
 &\phi''(y)= \frac{\sigma^2}{\eta^2} \left(\phi(y)\right)^2,\quad y>0,\\
&\phi(0)=1, \quad \lim_{y\to\infty}\phi(y)=0,
\end{array}\right.
$$
is given by $\left(\frac{\sigma}{\sqrt{6}\eta} y+1\right)^{-2}$, which
plays an important role and leads to the limit behavior \eqref{Tail-probability}. In our case, the above problem is replaced by  the following problem:
$$
\left\{\begin{array}{rl}
 &\phi''(y)= C \left(\phi(y)\right)^{\alpha},\quad y>0,\\
&\phi(0)=1, \quad \lim_{y\to\infty}\phi(y)=0
\end{array}\right.
$$
with $C$ being some positive constant. The solution to the above problem is $(\theta y+1)^{-\frac{2}{\alpha-1}}$ with some constant $\theta>0$ (see the proof of Corollary \ref{cor1}), which leads to the limit behavior \eqref{main}.

\section{Preliminaries}\label{pre}

Set $\widetilde{\xi}_t:= - \xi_t$.  Consider a branching L\'{e}vy process $\{ \widetilde{X}_u(t), u\in \widetilde{N}_t, t>0 \}$ with spatial motion $\widetilde{\xi}$, branching rate $\beta>0$ and  offspring distribution $\{p_k: k\geq 0\}$. Then
\begin{align}\label{step_1}
	\P \left(M < x\right) = \P \left( \inf_{t\geq 0} \inf_{u\in N_t} \widetilde{X}_u(t) >  -x\right) = \P_x \left( \inf_{t\geq 0} \inf_{u\in N_t} \widetilde{X}_u(t) > 0\right),
\end{align}
with the convention $\inf_{u\in N_t} \widetilde{X}_u(t) = +\infty$ when $\widetilde{N}_t=\emptyset$.
Define
\begin{align}
	v(x):= \P(M\geq x) \quad \mbox{and}\quad \widetilde{\tau}_y:= \inf\left\{t>0: \widetilde{\xi}_t \leq y\right\}.
\end{align}
It is easily seen that $v(x)=1$ for $x\leq 0$.

\subsection{Moment for overshoot of L\'{e}vy process}

For integer-valued random walks, the following result can be found in \cite[Lemma 10] {LS15}. We now prove that it also holds for some L\'{e}vy processes.

\begin{lemma}\label{lemma4}
	Let $\widetilde{\xi}_t$ be a L\'{e}vy process with $\Pi_0\left(\widetilde{\xi}_1\right)=0$ and $\Pi_0\left(\widetilde{\xi}_1^2\right)<\infty$. If $\Pi_0\left(|\widetilde{\xi}_1|^r\right)<\infty$ for some $r>2$, then
	\begin{align}
		\sup_{x>0} \Pi_x\left(\left|\widetilde{\xi}_{\widetilde{\tau}_0}\right|^{r-2} \right)<\infty.
	\end{align}
\end{lemma}
\textbf{Proof: }  Assume that $\Pi_0 (e^{\mathrm{i}\theta \widetilde{\xi}_1})= e^{-\Psi(\mathrm{i}\theta)}$ where
\[
\Psi (\mathrm{i}\theta ) =-\mathrm{i} \gamma \theta +\frac{\nu^2}{2}\theta^2+ \int_{x\neq 0}\left(1-e^{\mathrm{i}\theta x}+ \mathrm{i}\theta x1_{\{|x|\in (0,1]\}}\right)\pi(\mathrm{d}x)
\]
with $\pi$ being the L\'{e}vy measure.

(i) If $\pi(\{|x| >1\})=0$,  then by \cite[Theorem 36.7]{Sato}, $\widetilde{\xi}$ is recurrent and so
\[
\sup_{x>0} \Pi_x\left(\left|\widetilde{\xi}_{\widetilde{\tau}_0}\right|^{r-2} \right)\leq 1<\infty.
\]

(ii) If $\pi(\{ |x|>1\})>0$,  let $\sigma_n$ be the $n$-th time that $\widetilde{\xi}$ has a
jump of magnitude larger than $1$.  Similar to \cite[p.208]{DM02}, for $j\ge 1$, define $W_j=\widetilde{\xi}_{\sigma_j-}$ and
$V_j= \widetilde{\xi}_{\sigma_j}-\widetilde{\xi}_{\sigma_j-}$.
   Then $\{W_j: j\ge 1 \}$ and $\{V_j: j\ge 1\}$ are both iid families of random variables and independent of each other.  Furthermore,
\begin{align}\label{Distribution-Jump}
\Pi_0(V_1\in \mathrm{d}x) = \frac{\pi(\mathrm{d}x)}{ \pi(\{ |x|>1\})} 1_{\{|x|>1\}}
\end{align}
and
$W_1 \stackrel{\mathrm{d}}{=} \widetilde{\xi}^{(1)}_\mathbf{e}$ where $\widetilde{\xi}^{(1)}$ is a L\'{e}vy process with
\[
\Pi_0\left(e^{\mathrm{i}\theta \widetilde{\xi}^{(1)}_1} \right) =\exp\left\{ \mathrm{i} \gamma \theta -\frac{\nu^2}{2}\theta^2- \int_{|x|\in (0,1]}\left(1-e^{\mathrm{i}\theta x}+ \mathrm{i}\theta x1_{\{|x|\in (0,1]\}}\right)\pi(\mathrm{d}x)\right\}
\]
and $\mathbf{e}$ is an independent exponential random variable with parameter $\pi(\{ |x|>1\})$.
Therefore, by \eqref{Distribution-Jump} and  \cite[Theorem 25.3]{Sato},
\begin{align}\label{step_27}
	\Pi_0\left(|\widetilde{\xi}_1|^r\right)<\infty \quad \Longleftrightarrow \quad \int_{|x|>1} |x|^r \pi(\mathrm{d}x)<\infty \quad \Longleftrightarrow  \quad \Pi_0\left(|V_1|^r\right)<\infty.
\end{align}
By the definition of $W_1$, we know that  $\Pi_0\left(|W_1|^r\right)<\infty$. For $n\ge 1$, let
$$
Z_n:= \widetilde{\xi}_{\sigma_n}= \sum_{j=1}^n \left(W_j+ V_j\right).
$$
Then combining the above with \eqref{step_27}, we get
\[
	\Pi_0\left(|\widetilde{\xi}_1|^r\right)<\infty \quad \Longleftrightarrow \quad \Pi_0\left( |Z_1|^r \right)<\infty.
\]
By \cite[p.209]{DM02},
 for all $z>1$ and $x>0$,
\[
\Pi_x \left( |\widetilde{\xi}_{\widetilde{\tau}_0}| >z\right)\leq \Pi_x\left( |Z_{\widehat{\tau}_0}|>z\right),
\]
where $\widehat{\tau}_0:=\inf\{n:  Z_n<0\}$. Define
\begin{align}
	T_1:= \min\{n>0: Z_n< Z_0\},\quad T_k:=\inf\{n> T_{k-1}: Z_n< Z_{T_{k-1}}\},\quad S_n:= Z_{T_n},
\end{align}
then $S_1, S_2-S_1, S_3-S_2, \dots$,  are iid with finite  $(r-1)$-th moment if $\Pi_0\left( |Z_1|^r \right)<\infty$ (see \cite[Corollary 1]{Doney80}).
Note that for $z>1$,
\begin{align}\label{step_34}
	&\Pi_x\left(|Z_{\widehat{\tau}_0}|>z\right)= \sum_{k=0}^\infty \Pi_x\left(S_k>0, S_{k+1}<-z  \right)\nonumber\\
	&\leq \sum_{\ell=0}^{[x]} \left(\sum_{k=0}^\infty \Pi_x\left(S_k\in [\ell, \ell+1]\right) \right)	\Pi_0 \left(|S_1|> z+ \ell\right).
\end{align}
For any $\ell \in \N$, set $\tau^{(\ell)}:= \inf\{n: S_n \leq \ell+1\}$. Note that on the set $\left\{\sum_{k=0}^\infty 1_{\{S_k\in[ \ell,  \ell+1] \}} \geq m\right\}$, we have
    $S_{\tau^{(\ell)}+m-1}- S_{\tau^{(\ell)}}\in [-1,0]$.
Thus by the strong Markov property,
\begin{align}\label{e:rsnew}
     \Pi_x\left(\sum_{k=0}^\infty 1_{\{S_k\in[ \ell,  \ell+1] \}} \geq m\right) \leq \Pi_x\left( S_{\tau+m-1}- S_{\tau}\in [-1,0]\right)=
     \Pi_0\left( S_{m-1} \in [-1, 0]\right).
\end{align}
Consequently
\begin{align}\label{step_35}
	& \sum_{k=0}^\infty \Pi_x\left(S_k\in [\ell, \ell+1]\right) = \sum_{m=0}^\infty m \Pi_x\left( \sum_{k=0}^\infty 1_{\{S_k\in[ \ell,  \ell+1] \}} =m\right)\nonumber\\
	& \leq 1+ \sum_{m=0}^\infty \Pi_x\left(\sum_{k=0}^\infty 1_{\{S_k\in[ \ell,  \ell+1] \}} \geq m\right)
 	\leq 2+ \sum_{m=1}^\infty \Pi_0\left( S_{m-1} \in [-1, 0] \right)=: C<\infty,
\end{align}
where in the second to last inequality we sued \eqref{e:rsnew} and in the last we used the fact that $S_n\to -\infty$. Combining \eqref{step_34} and \eqref{step_35}, we get that
\begin{align}
&	\sup_{x>0} \Pi_x\left(|Z_{\widehat{\tau}_0}|>z\right) \leq C \sum_{\ell=0}^\infty \Pi_0 \left(|S_1|> z+ \ell\right)\nonumber\\
&\leq C\int_0^\infty\Pi_0 \left(|S_1|> z+ y -1\right)\mathrm{d} y \leq C \Pi_0\left(|S_1|1_{\{|S_1|> z-1\}}\right).
\end{align}
Therefore,
\begin{align}
& 	\Pi_x\left(\left|\widetilde{\xi}_{\widetilde{\tau}_0}\right|^{r-2} \right) \leq 2^{r-2} + (r-2) \int_2^\infty z^{r-3} \Pi_x \left( |\widetilde{\xi}_{\widetilde{\tau}_0}| >z\right)\mathrm{d} z\nonumber\\
	& \leq 2^{r-2}+ (r-2) \int_2^\infty z^{r-3} \Pi_x \left(|Z_{\widehat{\tau}_0}|>z\right)\mathrm{d} z \nonumber\\
	& \leq 2^{r-2}+ C(r-2) \int_2^\infty z^{r-3}   \Pi_0\left(|S_1|1_{\{|S_1|> z-1\}}\right)  \mathrm{d} z \leq 2^{r-2} + C \Pi_0\left(|S_1| (|S_1|+1)^{r-2}\right)<\infty,
\end{align}
which completes the proof of the lemma.
\hfill$\Box$

\subsection{Feynman-Kac representation for $v(x)$}

Define a function $f: [0, 1]\mapsto \R$ by
\begin{align}\label{Function-F}
	f(v):= \beta \frac{\sum_{k=0}^\infty p_k(1-v)^k - (1-v)}{v},\quad v\in (0,1],
\end{align}
and $f(0):= f(0+)=0$.
It is easy to see that $f(v)\geq 0$ for $v\in [0,1]$.
Also, define
\[
F(v)=\frac{1}{v}\left(1-\sum_{k=0}^\infty p_k (1-v)^k\right), \quad v\in (0,1].
\]
Note that $\beta(F(v)-1)=-f(v)$.
Recall that $v(x)= \P(M\geq x)$.

\begin{lemma}\label{Feynman-Kac}
	 For any $0\leq y<x$,
	 \begin{align}
	 	v(x)= \Pi_x\left(\exp\left\{ -\int_0^{\widetilde{\tau}_y} f\left(v\left( \widetilde{\xi}_s\right)\right)\mathrm{d}s\right\} v\left(\widetilde{\xi}_{\widetilde{\tau}_y}\right)\right).
	 \end{align}
\end{lemma}
\textbf{Proof: }
Put $u(x)=1-v(x)$.
Since the first branching time is  an independent
exponential random variable of parameter $\beta$, we have
\begin{align}
	& u(x)=  \P_x \left( \inf_{t\geq 0} \inf_{u\in N_t} \widetilde{X}_u(t) > 0\right)=\int_0^\infty \beta e^{-\beta s}\sum_{k=0}^\infty p_k \Pi_x \left(1_{\{\widetilde{\tau}_0 >s\}} \left(u(\widetilde{\xi}_s)\right)^k\right)\mathrm{d}s\nonumber\\
	& =\Pi_x\left( \int_0^{\widetilde{\tau}_0} \beta e^{-\beta s}\sum_{k=0}^\infty p_k  \left(u(\widetilde{\xi}_s)\right)^k\mathrm{d}s\right).
\end{align}
According to \cite[Lemma 4.1]{Dy01}, we have
\begin{align}
	u(x)+\beta \Pi_x \left(\int_0^{\widetilde{\tau}_0} u(\widetilde{\xi}_s)\mathrm{d}s\right)=\beta \Pi_x\left( \int_0^{\widetilde{\tau}_0}\sum_{k=0}^\infty p_k  \left(u(\widetilde{\xi}_s)\right)^k\mathrm{d}s\right),
\end{align}
which is  equivalent to
\begin{align}
	v(x)= 1- \beta \Pi_x \left(\int_0^{\widetilde{\tau}_0} \sum_{k=0}^\infty p_k  \left(1-v(\widetilde{\xi}_s)\right)^k - \left(1-v(\widetilde{\xi}_s)\right)\mathrm{d}s \right)= 1-\Pi_x \left(\int_0^{\widetilde{\tau}_0} f(v(\widetilde{\xi}_s))v(\widetilde{\xi}_s)\mathrm{d}s \right),
\end{align}
which can be written as
$$v(x)+\Pi_x \left(\int_0^{\widetilde{\tau}_0} f(v(\widetilde{\xi}_s))v(\widetilde{\xi}_s)\mathrm{d}s \right)=1.$$
Therefore, $v$ is a solution of the Schr\"{o}dinger equation: $v(x)+\Pi_x \left(\int_0^{\widetilde{\tau}_0} c(\widetilde{\xi}_s)v(\widetilde{\xi}_s)\mathrm{d}s \right)=1$  in $(0, \infty)$ with  $c(x):=f(v(x))\geq 0$. Then we have
\begin{align}
	v(x)= \Pi_x\left(\exp\left\{ -\int_0^{\widetilde{\tau}_0} f\left(v\left( \widetilde{\xi}_s\right)\right)\mathrm{d}s\right\} \right).
\end{align}
The desired result follows by conditioning on $\mathcal{F}_{\widetilde{\tau}_y}$ and applying  the strong Markov property of $\widetilde{\xi}$.
\hfill$\Box$

\subsection{An invariance principle for L\'{e}vy process}

The following lemma is an invariance principle for L\'{e}vy process. The proof is standard by comparing with random walks. We omit the proof here.

\begin{lemma}\label{lemma3}
	Suppose that $\widetilde{\xi}_t$ is a L\'{e}vy process with $\Pi_0 (\widetilde{\xi}_1)=0, \eta^2= \Pi_0(\widetilde{\xi}_1^2)\in (0, \infty)$.
	Assume that there exists $\varepsilon>0$ such that $\Pi_0(|\widetilde{\xi}_1|^{2+\varepsilon})<\infty$.
Then the processes
	\[
	\frac{\widetilde{\xi}_{nt}}{\eta \sqrt{n}},\quad t\in [0,\infty)
	\]
	converges weakly to a standard Brownian motion
$\{B_t, t\geq 0\}$
in the Skorohod topology.
\end{lemma}

\section{Proof of the main result}\label{Main}

\begin{lemma}\label{lemma2}
	Under the assumption \eqref{condition-1},  the function $f$ defined in \eqref{Function-F} satisfies that
	\begin{align}
		 \lim_{v\downarrow 0} \frac{f(v)}{v^{\alpha-1}}= \frac{\beta \kappa \Gamma(2-\alpha)}{\alpha-1}.
	\end{align}
\end{lemma}
\textbf{Proof: }
 Let $L$ be a random variable with the offspring distribution $\{p_k; k\ge 0\}$. It follows from  \cite[Theorem 8.1.6]{BGT}
 that $\mathbb{P}(L>x)\stackrel{x\to+\infty}{\sim}x^{-\alpha} c$ is equivalent to $\mathbb{E}(e^{-s L}) -1+ \mathbb{E}(L)s\stackrel{s\to 0}{\sim}s^\alpha \frac{\Gamma(2-\alpha)}{\alpha-1}c$, which is in turn equivalent to $\mathbb{E}(e^{-s L}) -e^{-s \mathbb{E}(L)}\stackrel{s\to 0}{\sim}s^\alpha \frac{\Gamma(2-\alpha)}{\alpha-1}c$.
 Therefore, letting $1-v=e^{-s}$,	 \eqref{condition-1} is equivalent to
\begin{align}
	\lim_{v\downarrow 0} \frac{vf(v)}{(-\ln (1-v))^{\alpha}}= \frac{\beta \kappa \Gamma(2-\alpha)}{\alpha-1} ,
\end{align}
which completes the proof of the lemma
since $\lim_{v\downarrow 0} \frac{v^{\alpha}}{(-\ln (1-v))^{\alpha}}=1.$
\hfill$\Box$

For any fixed $y\geq 0$,  the function
\[
 [0,\infty) \ni x\mapsto \frac{v\left(x+ y v(x)^{-\frac{\alpha-1}{2}}\right)}{v(x)}
\]
is bounded between $0$ and $1$.  Therefore, by  a diagonalization argument, we can find a subsequence
$\{x_k\in [0,\infty)\}$ with $\lim_{k\to\infty} x_k =+\infty$ such that for all $y\geq 0, y\in \Q$, the following limits exist:
\begin{align}\label{step_7}
	\phi(y):= \lim_{k\to\infty} \frac{v\left(x_k+ y v(x_k)^{-\frac{\alpha-1}{2}}\right)}{v(x_k)}.
\end{align}
Using the fact that $v(x)$ is decreasing, we see that $\phi(0)=1$ and $\phi(y)\in [0,1]$ for any $y\in  \Q\cap [0,\infty)$. Moreover,  for non-negative rational numbers $y_1<y_2$,
it holds that $\phi(y_1)\geq \phi(y_2)$. Therefore, for any $y\geq 0$, we can define
\begin{align}\label{step_8}
	\phi(y):= \sup_{z\in \Q, z\geq y} \phi(z) = \lim_{z\in \Q, z\downarrow y} \phi(y).
\end{align}

\begin{prop}\label{prop1}
The function $\phi(y)$ is a continuous decreasing function in $[0,\infty)$ and
\begin{align}\label{step_7a}
\phi(y)= \lim_{k\to\infty} \frac{v\left(x_k+ y v(x_k)^{-\frac{\alpha-1}{2}}\right)}{v(x_k)}, \quad \mbox{for all } y\ge 0.
\end{align}
Moreover, for any $K>0$,
 we have uniformly for $y\in [0,K]$,
\begin{align}\label{uniform-limit}
 \lim_{k\to\infty} \frac{v\left(x_k+ y v(x_k)^{-\frac{\alpha-1}{2}}\right)}{\phi(y)v(x_k)}=1.
\end{align}
\end{prop}
\textbf{Proof:}
Fix two non-negative rational numbers $y_1< y_2$.
By Lemma \ref{lemma2}, there exists a constant $C_1>0$ such that $f(v)\leq C_1 v^{\alpha-1}$ for all $v\in [0,1]$.
Set $z_i(k)= y_i v(x_k)^{-\frac{\alpha-1}{2}}$.  It follows from  Lemma \ref{Feynman-Kac} that
\begin{align}\label{step_9}
	&\phi(y_1)\geq  \phi(y_2)=  \lim_{k\to\infty} \frac{v\left(x_k+ z_2(k)\right)}{v(x_k)}\nonumber\\
	& =  \lim_{k\to\infty}  \Pi_{x_k+ z_2 (k)}\left(\exp\left\{ -\int_0^{\widetilde{\tau}_{x_k+ z_1(k)} } f\left(v\left( \widetilde{\xi}_s\right)\right)\mathrm{d}s\right\} \frac{v\left(\widetilde{\xi}_{\widetilde{\tau}_{x_k+z_1(k)}}\right) }{v(x_k)}\right)\nonumber\\
	&\geq \limsup_{k\to\infty}  \Pi_{x_k+ z_2 (k)}\left(\exp\left\{ - C_1 \int_0^{\widetilde{\tau}_{x_k+ z_1(k)} } \left(v\left( \widetilde{\xi}_s\right)\right)^{\alpha-1}\mathrm{d}s\right\} \right) \frac{v\left(x_k + z_1(k)\right) }{v(x_k)},
\end{align}
where in the last inequality, we used the fact that $v$ is decreasing and that $\widetilde{\xi}_{\widetilde{\tau}_{x_k+z_1(k)}}\leq x_k+z_1(k).$  Since $\widetilde{\xi}_s\geq x_k+z_1(k) \geq x_k$ for $s\in (0, \widetilde{\tau}_{x_k+ z_1(k)} ) $ and $v$ is decreasing, by \eqref{step_9}, we have
\begin{align}\label{step_10}
	&\phi(y_1)\geq  \phi(y_2)\geq \phi(y_1) \limsup_{k\to\infty}  \Pi_{x_k+ z_2 (k)}\left(\exp\left\{ - C_1 \left(v\left( x_k\right)\right)^{\alpha-1} \widetilde{\tau}_{x_k+ z_1(k)}  \right\} \right) \nonumber\\
	&=  \phi(y_1) \limsup_{k\to\infty}  \Pi_{0}\left(\exp\left\{ - C_1 \left(v\left( x_k\right)\right)^{\alpha-1} \widetilde{\tau}_{z_1(k)-z_2(k)}  \right\} \right)  .
\end{align}
Set $a:= y_2-y_1>0, n_k:= \left(v\left( x_k\right)\right)^{-(\alpha-1)}$. Since for $t>0$,
\begin{align}
	\Pi_0\left( n_k^{-1} \widetilde{\tau}_{-a n_k^{1/2}} >t\right) = \Pi_0\left( n_k^{-1/2}\inf_{s\leq tn_k} \widetilde{\xi}_s > -a \right) = \Pi_0\left(\frac{\inf_{s\leq t} \widetilde{\xi}_{n_k s} }{n_k^{1/2}}> -a \right),
\end{align}
it follows from Lemma \ref{lemma3} that
\begin{align}\label{step_11}
	&\lim_{k\to\infty}\Pi_0\left( n_k^{-1} \widetilde{\tau}_{-a n_k^{1/2}} >t\right)= \lim_{k\to\infty}\Pi_0\left(\inf_{s\leq  t} \frac{\widetilde{\xi}_{n_k s}}{n_k^{1/2}}> -a \right) \nonumber\\
	& = \Pi_0\left(\eta \inf_{s\leq t} B_t> -a \right) = \Pi_0\left(\tau_{-a\eta^{-1}}^{BM}>t \right),
\end{align}
where $\tau_b^{BM}$ is the first time that a standard Brownian motion hits $b$.
Combining \eqref{step_10} and \eqref{step_11},
\begin{align}\label{step_12}
	&\phi(y_1)\geq  \phi(y_2) \geq \phi(y_1) \Pi_{0}\left(\exp\left\{ - C_1 \tau_{\left(y_1-y_2\right)\eta^{-1}}^{BM}  \right\} \right) = e^{-\sqrt{2C_1} \frac{(y_2-y_1)}{\eta}}\phi(y_1)  .
\end{align}
By the definition of $\phi$ in \eqref{step_8}, we see that \eqref{step_12} holds
for all non-negative real numbers $y_1<y_2$.
This  implies that $\phi$ is continuous.
Besides, for any $y\geq 0$, we can fix two
non-negative rational numbers $y_1\leq y<y_2$.
 Then by the monotonicity of $v$,
\begin{align}
	& \phi(y_2)= \lim_{k\to\infty} \frac{v\left(x_k+ y_2 v(x_k)^{-\frac{\alpha-1}{2}}\right)}{v(x_k)} \leq \liminf_{k\to\infty}\frac{v\left(x_k+ y  v(x_k)^{-\frac{\alpha-1}{2}}\right)}{v(x_k)} \nonumber\\
	& \leq \limsup_{k\to\infty}\frac{v\left(x_k+ y  v(x_k)^{-\frac{\alpha-1}{2}}\right)}{v(x_k)} \leq \lim_{k\to\infty} \frac{v\left(x_k+ y_1 v(x_k)^{-\frac{\alpha-1}{2}}\right)}{v(x_k)} = \phi(y_1),
\end{align}
which implies \eqref{step_7a} by letting $y_1\uparrow y$ and $y_2\downarrow y$.

Finally we prove the uniform convergence.  For any $\epsilon>0$, we can  find $y_0=0 <y_1<...<y_m=K$ such that
\begin{align}
	\sup_{1\leq i\leq m} \left| \phi(y_i)-\phi(y_{i-1}) \right| < \frac{\epsilon}{2}.
\end{align}
Now we can find a common $N$ such that for all $0\leq i\leq m$, when $k>N$,
\[
\left|\frac{v\left(x_k+ y_i v(x_k)^{-\frac{\alpha-1}{2}}\right)}{v(x_k)} -\phi(y_i)\right|<\frac{\epsilon}{2}.
\]
Therefore, for any
$i=1,\dots,m$
and $y\in [y_{i-1}, y_i]$, when $k>N$,
\begin{align}\label{step_13}
	& \phi(y)-\epsilon \leq \phi(y_{i-1})-\epsilon< \phi(y_i)-\frac{\epsilon}{2} <\frac{v\left(x_k+ y_i v(x_k)^{-\frac{\alpha-1}{2}}\right)}{v(x_k)}\leq \frac{v\left(x_k+ y v(x_k)^{-\frac{\alpha-1}{2}}\right)}{v(x_k)}\nonumber\\
	& \leq \frac{v\left(x_k+ y_{i-1}v(x_k)^{-\frac{\alpha-1}{2}}\right)}{v(x_k)} < \frac{\epsilon}{2} + \phi(y_{i-1})< \epsilon +\phi(y_i)\leq \epsilon +\phi(y).
\end{align}
Noticing that  $\phi(0)=1$ and $\phi(K)>0$ which holds by \eqref{step_12} with $y_1=0, y_2=K$, by \eqref{step_13}, we obtain the desired result \eqref{uniform-limit}.
\hfill$\Box$

Given Lemma \ref{lemma3}  and  Proposition \ref{prop1}, the following result seems trivial,
%YX but we will give a formal proof.
%RS but we will give an outline of proof.
but we will give a proof.
Recall that $n_k= v(x_k)^{-(\alpha-1)} $
 and $\eta= \sqrt{\Pi_0(\widetilde{\xi}_1^2)}$.

\begin{lemma}\label{lemma1}
 	For any $\theta>0, y>0$ and $z\geq y$,  it holds that
	\begin{align}\label{step_40}
		&\lim_{k\to\infty} \Pi_{0}\left(\exp\left\{ - \theta \int_0^{n_k^{-1}
\widetilde{\tau}_{-y\sqrt{n_k}}}
 \left(\frac{v\left(  \left( n_k^{-1/2}\widetilde{\xi}_{n_k s}+z\right)v(x_k)^{-\frac{\alpha-1}{2}} + x_k \right)}{v(x_k)}\right)^{\alpha-1}\mathrm{d}s\right\} \right)\nonumber\\
		& = \Pi_{0}\left(\exp\left\{ -\theta
\int_0^{\tau_{-y/\eta}^{BM}}
\left(\phi(\eta B_s+z)\right)^{\alpha-1}\mathrm{d}s\right\}\right),
	\end{align}
where $\tau_{-y/\eta}^{BM}$ is the first time that a standard Brownian motion hits $-y/\eta$.
\end{lemma}
\textbf{Proof:} For simplicity, we set
\[
\widetilde{\tau}^{(k)}:= n_k^{-1}
\widetilde{\tau}_{-y\sqrt{n_k}},
\quad \widetilde{\xi}^{(k)}_s:= \frac{\widetilde{\xi}_{n_ks}}{\sqrt{n_k}}.
\]
\textbf{Step 1}:
In this step, we  prove that for any $T, A>0$,
\begin{align}\label{step_36}
	&\lim_{k\to\infty} \Pi_{0}\left(\exp\left\{ - \theta \int_0^{\widetilde{\tau}^{(k)} \land T} \left(\phi\left( \widetilde{\xi}^{(k)}_s +z \right)\right)^{\alpha-1}\mathrm{d}s\right\} 1_{\{\sup_{s\in[0,T]} \widetilde{\xi}_s^{(k)} < A \}}\right)\nonumber\\
	& = \Pi_{0}\left(\exp\left\{ -\theta \int_0^{\tau_{-y\eta^{-1}}^{BM} \land T} \left(\phi(\eta B_s+z)\right)^{\alpha-1}\mathrm{d}s\right\}1_{\{ \eta \sup_{s\in[0,T]}B_s < A \}}\right).
\end{align}
For any interger $N>1$, define $t_i:= Ti/N, 1\leq i\leq N$.  Since $\phi$ is decreasing, it holds that
\begin{align}\label{step_42}
	& \int_0^{\widetilde{\tau}^{(k)} \land T} \left(\phi\left( \widetilde{\xi}^{(k)}_s +z \right)\right)^{\alpha-1}\mathrm{d}s = \sum_{i=1}^N \int_{t_{i-1}}^{t_i} \left(\phi\left( \widetilde{\xi}^{(k)}_s +z \right)\right)^{\alpha-1} 1_{\{s<\widetilde{\tau}^{(k)}  \}}\mathrm{d}s\nonumber\\
	& \geq  \sum_{i=1}^N \int_{t_{i-1}}^{t_i} \left(\phi\left( \sup_{s\in [t_{i-1}, t_i]} \widetilde{\xi}^{(k)}_s +z \right)\right)^{\alpha-1} 1_{\{t_i<\widetilde{\tau}^{(k)}  \}}\mathrm{d}s \nonumber\\
	&= \frac{T}{N}\sum_{i=1}^N  \left(\phi\left( \sup_{s\in [t_{i-1}, t_i]} \widetilde{\xi}^{(k)}_s +z \right)\right)^{\alpha-1} 1_{\{t_i<\widetilde{\tau}^{(k)} \}}.
\end{align}
It is easy to check that $\{t_i<\widetilde{\tau}^{(k)} \}=\{ \inf_{s\leq t_i}\widetilde{\xi}_s^{(k)} > -y\}.$  Also,
observe that the  functionals
\[
w\in D[0, T] \mapsto \sup_{s\in [t_{j-1},t_j]} w(s)\in \R, \quad i=1, \dots, N.
\]
are continuous with respect to the Skorohod topology.  Therefore, taking two sequences of continuous functions $h_\ell(x) \uparrow 1_{(-y,+\infty)}(x)$ and $j_\ell(x)\downarrow 1_{(-\infty, A)}(x)$, by Lemma \ref{lemma3} and \eqref{step_42}, we get that
\begin{align}
		&\limsup_{k\to\infty} \Pi_{0}\left(\exp\left\{ - \theta  \frac{T}{N}\sum_{i=1}^N  \left(\phi\left( \sup_{s\in [t_{i-1}, t_i]} \widetilde{\xi}^{(k)}_s +z \right)\right)^{\alpha-1} 1_{\{t_i<\widetilde{\tau}^{(k)} \}}\right\} 1_{\{\sup_{s\in[0,T]} \widetilde{\xi}_s^{(k)} < A \}} \right)\nonumber\\
		& \leq \limsup_{k\to\infty} \Pi_{0}\left(\exp\left\{ - \theta  \frac{T}{N}\sum_{i=1}^N  \left(\phi\left( \sup_{s\in [t_{i-1}, t_i]} \widetilde{\xi}^{(k)}_s +z \right)\right)^{\alpha-1} h_\ell \left(\inf_{s\leq t_i}\widetilde{\xi}_s^{(k)}\right)\right\} j_\ell\left(\sup_{s\in[0,T]} \widetilde{\xi}_s^{(k)} \right)\right)\nonumber\\
		& = \Pi_{0}\left(\exp\left\{ -\theta \frac{T}{N}\sum_{i=1}^N  \left(\phi\left( \eta \sup_{s\in [t_{i-1}, t_i]} B_s +z \right)\right)^{\alpha-1} h_\ell\left( \eta \inf_{s\leq t_i} B_s\right) \right\} j_\ell\left(\eta \sup_{s\in[0,T]}B_s\right)\right).
\end{align}
Letting $\ell \to +\infty$, we get
\begin{align}\label{step_37}
	&\limsup_{k\to\infty} \Pi_{0}\left(\exp\left\{ - \theta \int_0^{\widetilde{\tau}^{(k)} \land T} \left(\phi\left( \widetilde{\xi}^{(k)}_s +z \right)\right)^{\alpha-1}\mathrm{d}s\right\} 1_{\{\sup_{s\in[0,T]} \widetilde{\xi}_s^{(k)} < A \}} \right)\nonumber\\
	&\leq  \Pi_{0}\left(\exp\left\{ -\theta \frac{T}{N}\sum_{i=1}^N  \left(\phi\left( \eta \sup_{s\in [t_{i-1}, t_i]} B_s +z \right)\right)^{\alpha-1} 1_{\left\{t_i<\tau_{-y\eta^{-1}}^{BM}  \right\}} \right\}1_{\left\{ \eta \sup_{s\in[0,T]}B_s < A \right\}}\right).
\end{align}
Letting $N\to+\infty$ in \eqref{step_37}, we get
\begin{align}
	&\limsup_{k\to\infty} \Pi_{0}\left(\exp\left\{ - \theta \int_0^{\widetilde{\tau}^{(k)} \land T} \left(\phi\left( \widetilde{\xi}^{(k)}_s +z \right)\right)^{\alpha-1}\mathrm{d}s\right\} 1_{\{\sup_{s\in[0,T]} \widetilde{\xi}_s^{(k)} < A \}} \right)\nonumber\\
	&\leq  \Pi_{0}\left(\exp\left\{ -\theta
 \int_0^{\tau_{-y/\eta}^{BM} \land T}
  \left(\phi(\eta B_s+z)\right)^{\alpha-1}\mathrm{d}s\right\}1_{\{ \eta \sup_{s\in[0,T]}B_s < A \}}\right).
\end{align}
Using a similar argument, we can get
\begin{align}
	&\liminf_{k\to\infty} \Pi_{0}\left(\exp\left\{ - \theta \int_0^{\widetilde{\tau}^{(k)} \land T} \left(\phi\left( \widetilde{\xi}^{(k)}_s +z \right)\right)^{\alpha-1}\mathrm{d}s\right\} 1_{\{\sup_{s\in[0,T]} \widetilde{\xi}_s^{(k)} < A \}} \right)\nonumber\\
	&\ge  \Pi_{0}\left(\exp\left\{ -\theta
\int_0^{\tau_{-y/\eta}^{BM} \land T}
 \left(\phi(\eta B_s+z)\right)^{\alpha-1}\mathrm{d}s\right\}1_{\{ \eta \sup_{s\in[0,T]}B_s < A \}}\right).
\end{align}
Combining the two displays above, we get the desired conclusion of this step.

\textbf{Step 2}:  In this step, we prove that for any $T, A>0$,
\begin{align}\label{step_39}
		&\lim_{k\to\infty} \Pi_{0}\left(\exp\left\{ - \theta \int_0^{\widetilde{\tau}^{(k)} \land T} \left(\phi^{(k)}\left( \widetilde{\xi}^{(k)}_s +z \right)\right)^{\alpha-1}\mathrm{d}s\right\} 1_{\{\sup_{s\in[0,T]} \widetilde{\xi}_s^{(k)} < A \}}\right)\nonumber\\
	& = \Pi_{0}\left(\exp\left\{ -\theta
\int_0^{\tau_{-y/\eta}^{BM} \land T}
 \left(\phi(\eta B_s+z)\right)^{\alpha-1}\mathrm{d}s\right\}1_{\{ \eta \sup_{s\in[0,T]}B_s < A \}}\right),
\end{align}
where
\[
\phi^{(k)}(z):= \frac{v\left(  \left( z\right)v(x_k)^{-\frac{\alpha-1}{2}} + x_k \right)}{v(x_k)}.
\]
Note that on set $\{\sup_{s\in[0,T]} \widetilde{\xi}_s^{(k)} < A \}$,  for any $s< \widetilde{\tau}^{(k)} \land T$,
    it holds that   $\widetilde{\xi}_s^{(k)}+z\in (z-y, A+z)\subset [0,A+z].$
It follows from Proposition \ref{prop1} that, for any $\varepsilon>0$, there exists $K$ such that for any $k>K$ and $s\in\widetilde{\tau}^{(k)} \land T$,
\begin{align}
	(1-\varepsilon)\left(\phi\left( \widetilde{\xi}^{(k)}_s +z \right)\right)^{\alpha-1} \leq  \left(\phi^{(k)}\left( \widetilde{\xi}^{(k)}_s +z \right)\right)^{\alpha-1} \leq 	(1+\varepsilon)\left(\phi\left( \widetilde{\xi}^{(k)}_s +z \right)\right)^{\alpha-1}.
\end{align}
Therefore, by \eqref{step_36},
\begin{align}\label{step_38}
	&\limsup_{k\to\infty} \Pi_{0}\left(\exp\left\{ - \theta \int_0^{\widetilde{\tau}^{(k)} \land T} \left(\phi^{(k)}\left( \widetilde{\xi}^{(k)}_s +z \right)\right)^{\alpha-1}\mathrm{d}s\right\} 1_{\{\sup_{s\in[0,T]} \widetilde{\xi}_s^{(k)} < A \}}\right)\nonumber\\
	&\leq \lim_{k\to\infty} \Pi_{0}\left(\exp\left\{ - \theta (1-\varepsilon)\int_0^{\widetilde{\tau}^{(k)} \land T} \left(\phi\left( \widetilde{\xi}^{(k)}_s +z \right)\right)^{\alpha-1}\mathrm{d}s\right\} 1_{\{\sup_{s\in[0,T]} \widetilde{\xi}_s^{(k)} < A \}}\right)\nonumber\\
	& =  \Pi_{0}\left(\exp\left\{ -\theta (1-\varepsilon)
\int_0^{\tau_{-y/\eta}^{BM} \land T}
\left(\phi(\eta B_s+z)\right)^{\alpha-1}\mathrm{d}s\right\}1_{\{ \eta \sup_{s\in[0,T]}B_s < A \}}\right).
\end{align}
Letting $\varepsilon \downarrow 0$, we get
\begin{align}
		&\limsup_{k\to\infty} \Pi_{0}\left(\exp\left\{ - \theta \int_0^{\widetilde{\tau}^{(k)} \land T} \left(\phi^{(k)}\left( \widetilde{\xi}^{(k)}_s +z \right)\right)^{\alpha-1}\mathrm{d}s\right\} 1_{\{\sup_{s\in[0,T]} \widetilde{\xi}_s^{(k)} < A \}}\right)\nonumber\\
	& \leq  \Pi_{0}\left(\exp\left\{ -\theta
\int_0^{\tau_{-y/\eta}^{BM} \land T}
\left(\phi(\eta B_s+z)\right)^{\alpha-1}\mathrm{d}s\right\}1_{\{ \eta \sup_{s\in[0,T]}B_s < A \}}\right).
\end{align}
Using a similar argument, we can get
\begin{align}
		&\liminf_{k\to\infty} \Pi_{0}\left(\exp\left\{ - \theta \int_0^{\widetilde{\tau}^{(k)} \land T} \left(\phi^{(k)}\left( \widetilde{\xi}^{(k)}_s +z \right)\right)^{\alpha-1}\mathrm{d}s\right\} 1_{\{\sup_{s\in[0,T]} \widetilde{\xi}_s^{(k)} < A \}}\right)\nonumber\\
	& \ge  \Pi_{0}\left(\exp\left\{ -\theta
\int_0^{\tau_{-y/\eta}^{BM} \land T}
 \left(\phi(\eta B_s+z)\right)^{\alpha-1}\mathrm{d}s\right\}1_{\{ \eta \sup_{s\in[0,T]}B_s < A \}}\right).
\end{align}
Combining the two displays above, we get the desired conclusion of this step.

\textbf{Step 3}: In this step, we prove \eqref{step_40}. Noting that
\begin{align}
	& \lim_{T\to\infty} \lim_{A\to\infty} \Pi_{0}\left(\exp\left\{ -\theta
    \int_0^{\tau_{-y/\eta}^{BM} \land T}
 \left(\phi(\eta B_s+z)\right)^{\alpha-1}\mathrm{d}s\right\}1_{\{ \eta \sup_{s\in[0,T]}B_s < A \}}\right)\nonumber\\
	& =  \Pi_{0}\left(\exp\left\{ -\theta
\int_0^{\tau_{-y/\eta}^{BM}}
\left(\phi(\eta B_s+z)\right)^{\alpha-1}\mathrm{d}s\right\}\right),
\end{align}
it suffices to prove that
\begin{align}\label{step_41}
	&\lim_{T\to\infty} \limsup_{A\to\infty} \limsup_{k\to\infty} \bigg|\Pi_{0}\left(\exp\left\{ - \theta \int_0^{\widetilde{\tau}^{(k)} \land T} \left(\phi^{(k)}\left( \widetilde{\xi}^{(k)}_s +z \right)\right)^{\alpha-1}\mathrm{d}s\right\} 1_{\{\sup_{s\in[0,T]} \widetilde{\xi}_s^{(k)} <A \}}\right)\nonumber\\
	& \quad -\Pi_{0}\left(\exp\left\{ - \theta \int_0^{\widetilde{\tau}^{(k)} } \left(\phi^{(k)}\left( \widetilde{\xi}^{(k)}_s +z \right)\right)^{\alpha-1}\mathrm{d}s\right\} \right) \bigg|=0.
\end{align}
The proof for \eqref{step_41} is standard so we omit the details here. This implies the desired result.
\hfill$\Box$

\begin{prop}\label{prop2}
	The function $\phi$ defined in \eqref{step_7} satisfies the equation
	\begin{align}
		\phi(y)= \Pi_0 \left(\exp\left\{ -  \frac{\beta \kappa \Gamma(2-\alpha)}{\alpha-1}
\int_0^{\tau_{-y/\eta}^{BM}}
\left(\phi(\eta B_s +y)\right)^{\alpha-1}\mathrm{d}s \right\} \right),\quad y\geq 0.
	\end{align}
\end{prop}
\textbf{Proof: }  Fix a constant $\rho>0$ and set $z_k:=x_k+v(x_k)^{-\frac{\alpha-1}2+\rho}$.  For $y>0$,  by Lemma \ref{Feynman-Kac},  we have
\begin{align}\label{step_14}
	& \frac{v(x_k+yv(x_k)^{-\frac{\alpha-1}{2}} +v(x_k)^{-\frac{\alpha-1}{2}+\rho} )}{v(x_k)} = \frac{v(z_k + yv(x_k)^{-\frac{\alpha-1}{2}})}{v(x_k)} \nonumber\\
	& =  \Pi_{z_k+yv(x_k)^{-\frac{\alpha-1}{2}} }\left(\exp\left\{ -\int_0^{\widetilde{\tau}_{z_k}} f\left(v\left( \widetilde{\xi}_s\right)\right)\mathrm{d}s\right\} \frac{v\left(\widetilde{\xi}_{\widetilde{\tau}_{z_k}}\right)}{v(x_k)}\right)\nonumber\\
	& = \Pi_{yv(x_k)^{-\frac{\alpha-1}{2}} }\left(\exp\left\{ -\int_0^{\widetilde{\tau}_0} f\left(v\left( \widetilde{\xi}_s + z_k\right)\right)\mathrm{d}s\right\} \frac{v\left(\widetilde{\xi}_{\widetilde{\tau}_0}+z_k\right)}{v(x_k)}\right).
\end{align}

We first show that
\begin{align}\label{step_15}
	\lim_{k\to\infty} \Pi_{yv(x_k)^{-\frac{\alpha-1}{2}} }\left( \left| \frac{v\left(\widetilde{\xi}_{\widetilde{\tau}_0}+z_k\right)}{v(x_k)} -1\right|\right)=0.
\end{align}
Indeed, on the event
\[
A:= \left\{\widetilde{\xi}_{\widetilde{\tau}_0} +z_k \geq x_k \right\},
\]
by the inequality $v(x_k)\geq v\left(\widetilde{\xi}_{\widetilde{\tau}_0}+z_k\right)\geq v(z_k)$, we have
\[
\left| \frac{v\left(\widetilde{\xi}_{\widetilde{\tau}_0}+z_k\right)}{v(x_k)} -1\right| = 1- \frac{v\left(\widetilde{\xi}_{\widetilde{\tau}_0}+z_k\right)}{v(x_k)} \leq 1- \frac{v\left(z_k\right)}{v(x_k)} ,
\]
 and on $A^c$, we have
 \[
 \left| \frac{v\left(\widetilde{\xi}_{\widetilde{\tau}_0}+z_k\right)}{v(x_k)} -1\right| \leq \frac{2}{v(x_k)}.
 \]
Therefore,
\begin{align}\label{step_17}
	\Pi_{yv(x_k)^{-\frac{\alpha-1}{2}} }\left( \left| \frac{v\left(\widetilde{\xi}_{\widetilde{\tau}_0}+z_k\right)}{v(x_k)} -1\right|\right) \leq \frac{2}{v(x_k)} \Pi_{yv(x_k)^{-\frac{\alpha-1}{2}} }\left(A^c\right)+1 - \frac{v\left(z_k\right)}{v(x_k)} .
\end{align}
By Markov's inequality, for any $r>2$, we have
\begin{align}
	\frac{1}{v(x_k)} \Pi_{yv(x_k)^{-\frac{\alpha-1}{2}} }\left(A^c\right) \leq  \Pi_{yv(x_k)^{-\frac{\alpha-1}{2}} } \left(\left| \widetilde{\xi}_{\widetilde{\tau}_0}\right|^{r-2}\right) \left(v(x_k)\right)^{\left(\frac{\alpha-1}{2}-\rho\right)(r-2)-1}.
\end{align}
Since $r>2\alpha/(\alpha-1)$, we can find a sufficient small $\rho>0$ such that $\left(\frac{\alpha-1}{2}-\rho\right)(r-2)>1$. Therefore, by Lemma \ref{lemma4}, we have
\begin{align}\label{step_18}
	\lim_{k\to\infty} 	\frac{1}{v(x_k)} \Pi_{yv(x_k)^{-\frac{\alpha-1}{2}} }\left(A^c\right)  =0.
\end{align}
Since $\lim_{k\to\infty} v(z_k)/v(x_k) = 1$ by Proposition \ref{prop1}, we immediately get \eqref{step_15} by combining \eqref{step_17} and \eqref{step_18}.

Letting $k\to\infty$, the left-hand side of \eqref{step_14} converges to $\phi(y)$ according to Proposition \ref{prop1}.  For the right-hand side of \eqref{step_14}, combining \eqref{step_15} and the trivial inequality $\left|\E(e^{-|X|} Y)- \E (e^{-|X|})\right| \leq \E \left(|Y -1|\right)$, we get that
\begin{align}\label{step_16}
	&\phi(y)= \lim_{k\to\infty} \Pi_{yv(x_k)^{-\frac{\alpha-1}{2}} }\left(\exp\left\{ -\int_0^{\widetilde{\tau}_0} f\left(v\left( \widetilde{\xi}_s + z_k\right)\right)\mathrm{d}s\right\} \right).
\end{align}
Using Lemma \ref{lemma2} and the fact that $\sup_{s<\widetilde{\tau}_0} v\left( \widetilde{\xi}_s + z_k\right)\leq v(z_k)\to 0$, we get that
for any $\varepsilon>0$, there exists $N$ such that for all $k\geq N$ and $s\in (0,  \widetilde{\tau}_0)$,
\begin{align}
 & \frac{\beta \kappa \Gamma(2-\alpha)}{\alpha-1}(1-\varepsilon)  \left(v\left( \widetilde{\xi}_s + x_k +  \varepsilon v(x_k)^{-\frac{\alpha-1}{2}}\right)\right)^{\alpha-1} \leq f\left(v\left( \widetilde{\xi}_s + z_k\right)\right) \nonumber\\
 &\leq  \frac{\beta \kappa \Gamma(2-\alpha)}{\alpha-1}(1+\varepsilon)  \left(v\left( \widetilde{\xi}_s + x_k \right)\right)^{\alpha-1} .
\end{align}
Plugging this into \eqref{step_16}, we get that
\begin{align}
	&\phi(y)\nonumber\\
	&\leq \liminf_{k\to\infty} \Pi_{yv(x_k)^{-\frac{\alpha-1}{2}} }\left(\exp\left\{ - \frac{\beta \kappa \Gamma(2-\alpha)}{\alpha-1}(1-\varepsilon)  \int_0^{\widetilde{\tau}_0} \left(v\left( \widetilde{\xi}_s + x_k + \varepsilon v(x_k)^{-\frac{\alpha-1}{2}}\right)\right)^{\alpha-1}\mathrm{d}s\right\} \right).
\end{align}
Note that for $n_k= v(x_k)^{-(\alpha-1)} $,
\begin{align}
	& \Pi_{yv(x_k)^{-\frac{\alpha-1}{2}} }\left(\exp\left\{ - \frac{\beta \kappa \Gamma(2-\alpha)}{\alpha-1}(1-\varepsilon)  \int_0^{\widetilde{\tau}_0} \left(v\left( \widetilde{\xi}_s + x_k + \varepsilon v(x_k)^{-\frac{\alpha-1}{2}}\right)\right)^{\alpha-1}\mathrm{d}s\right\} \right)\nonumber\\
	& = \Pi_{0}\left(\exp\left\{ - \frac{\beta \kappa \Gamma(2-\alpha)}{\alpha-1}(1-\varepsilon)  \int_0^{\widetilde{\tau}_{-yn_k^{1/2}}} \left(v\left( \widetilde{\xi}_s + (y+\varepsilon)v(x_k)^{-\frac{\alpha-1}{2}} + x_k \right)\right)^{\alpha-1}\mathrm{d}s\right\} \right)\nonumber\\
	& =  \Pi_{0}\left(\exp\left\{ - \frac{\beta \kappa \Gamma(2-\alpha)}{\alpha-1}(1-\varepsilon)  \int_0^{n_k^{-1}\widetilde{\tau}_{-yn_k^{1/2}}} \left(\frac{v\left(  \left( n_k^{-1/2}\widetilde{\xi}_{n_k s}+y+\varepsilon\right)v(x_k)^{-\frac{\alpha-1}{2}} + x_k \right)}{v(x_k)}\right)^{\alpha-1}\mathrm{d}s\right\} \right).
\end{align}
By Lemma \ref{lemma1},
\begin{align}
	& \lim_{k\to\infty} \Pi_{yv(x_k)^{-\frac{\alpha-1}{2}} }\left(\exp\left\{ - \frac{\beta \kappa \Gamma(2-\alpha)}{\alpha-1}(1-\varepsilon)  \int_0^{\widetilde{\tau}_0} \left(v\left( \widetilde{\xi}_s + x_k + \varepsilon v(x_k)^{-\frac{\alpha-1}{2}}\right)\right)^{\alpha-1}\mathrm{d}s\right\} \right)\nonumber\\
	& = \Pi_{0}\left(\exp\left\{ - \frac{\beta \kappa \Gamma(2-\alpha)}{\alpha-1}(1-\varepsilon)  \int_0^{\tau_{-y/\eta}^{BM}} \left(\phi(\eta B_s+y+\varepsilon)\right)^{\alpha-1}\mathrm{d}s\right\}\right).
\end{align}
Therefore, we conclude that
\begin{align}
		&\phi(y)\leq \Pi_{0}\left(\exp\left\{ - \frac{\beta \kappa \Gamma(2-\alpha)}{\alpha-1}(1-\varepsilon)  \int_0^{\tau_{-y\eta^{-1}}^{BM}} \left(\phi(\eta B_s+y+\varepsilon)\right)^{\alpha-1}\mathrm{d}s\right\}\right).
\end{align}
Let $\varepsilon\downarrow 0$, we obtain that
\begin{align}
	&\phi(y)\leq \Pi_{0}\left(\exp\left\{ - \frac{\beta \kappa \Gamma(2-\alpha)}{\alpha-1} \int_0^{\tau_{-y/\eta}^{BM}} \left(\phi(\eta B_s+y)\right)^{\alpha-1}\mathrm{d}s\right\}\right).
\end{align}
Similarly, we also have
\begin{align}
	&\phi(y)\geq \Pi_{0}\left(\exp\left\{ - \frac{\beta \kappa \Gamma(2-\alpha)}{\alpha-1} \int_0^{\tau_{-y/\eta}^{BM}} \left(\phi(\eta B_s+y)\right)^{\alpha-1}\mathrm{d}s\right\}\right).
\end{align}
Combining the two displays above, we arrive at the desired result.
\hfill$\Box$

\begin{cor}\label{cor1}
	It holds that
	\[
	\phi(y) =(\theta y+1)^{-\frac{2}{\alpha-1}},
	\]
	where
	\[
	\theta := \left( \frac{\beta\kappa \Gamma(2-\alpha)(\alpha-1)}{\eta^2 (\alpha+1)}\right)^{1/2}.
	\]
\end{cor}
\textbf{Proof: } By Proposition \ref{prop2}, $\phi$ is the unique solution to
$$
\left\{\begin{array}{rl}
 &\frac{\eta^2}{2}\phi''(y)= \frac{\beta \kappa \Gamma(2-\alpha)}{\alpha-1} \left(\phi(y)\right)^{\alpha},\quad y>0.\\
&\phi(0)=1, \quad \lim_{y\to\infty}\phi(y)=0.
\end{array}\right.
$$
It is easy to check that $\phi(y)=(\theta y+1)^{-\frac{2}{\alpha-1}}$  solves the above equation.

\hfill$\Box$

\textbf{Proof of Theorem \ref{thm1}} By Corollary \ref{cor1},  the limit $\phi$ is independent of  $\{x_k\}$, which implies that for all $y\geq 0$,
\begin{align}\label{step_19}
		(\theta y+1)^{-\frac{2}{\alpha-1}}= \lim_{x\to +\infty} \frac{v\left(x+ y v(x)^{-\frac{\alpha-1}{2}}\right)}{v(x)}.
\end{align}
Set $w(x)= x^{2/(\alpha-1)} v(x)$. Then \eqref{step_19} is equivalent to
\begin{align}\label{step_20}
	\lim_{x\to +\infty} \frac{w\left(x\left(1+ y w(x)^{-\frac{\alpha-1}{2}}\right)\right)}{w(x)} \cdot
	\frac{(\theta y+1)^{2/(\alpha-1)} }{\left( 1+ yw(x)^{-\frac{\alpha-1}{2}}\right)^{2/(\alpha-1)}}=1.
\end{align}
Suppose that
\begin{align}
	0\leq A := \liminf_{x\to \infty} w(x) \leq  \limsup_{x\to \infty} w(x) =: B\leq \infty.
\end{align}

\textbf{Step 1}: In this step, we prove $B>0$ and $A<\infty$.  Assume that $B=0$. In this case, for $k \in \N$,  define $b_k:= \sup\{x: w(x)> k^{-1}\}$, then $b_k\to +\infty$ and $w(b_k)=k^{-1}$. Taking $x=b_k$ and $y=1$ in \eqref{step_20}, we obtain that
\begin{align}
	\lim_{k\to +\infty} \frac{w\left(b_k\left(1+ k^{\frac{\alpha-1}{2}}\right)\right)}{k^{-1}} \cdot
	\frac{(\theta +1)^{2/(\alpha-1)} }{\left( 1+ k^{\frac{\alpha-1}{2}}\right)^{2/(\alpha-1)}}=1.
\end{align}
However, by the definition of $b_k$,
\begin{align}
	 \frac{w\left(b_k\left(1+ k^{\frac{\alpha-1}{2}}\right)\right)}{k^{-1}} \cdot
	\frac{(\theta +1)^{2/(\alpha-1)} }{\left( 1+ k^{\frac{\alpha-1}{2}}\right)^{2/(\alpha-1)}}\leq \frac{(\theta +1)^{2/(\alpha-1)} }{\left( 1+ k^{\frac{\alpha-1}{2}}\right)^{2/(\alpha-1)}} \stackrel{k\to\infty}{\longrightarrow}0,
\end{align}
which is a contradiction.  The proof of $A<\infty$ is similar.

\textbf{Step 2}:  In this step, we prove $A\leq \theta^{-2/(\alpha-1)}\leq B$.  By the definition of $B$, there exists $b_k \to+\infty$ such that $w(b_k)\to B$. Taking $x=b_k$ and $y=1$ in \eqref{step_20},  we get that
\begin{align}\label{step_21}
	\lim_{k\to +\infty} \frac{w\left(b_k\left(1+ k^{\frac{\alpha-1}{2}}\right)\right)}{B} \cdot
	\frac{(\theta +1)^{2/(\alpha-1)} }{\left( 1+ B^{-\frac{\alpha-1}{2}}\right)^{2/(\alpha-1)}}=1.
\end{align}
Since $\limsup_{k\to\infty} w\left(b_k\left(1+ k^{\frac{\alpha-1}{2}}\right)\right)\leq B$, \eqref{step_21} implies that
\[
1 \leq \frac{(\theta +1)^{2/(\alpha-1)} }{\left( 1+ B^{-\frac{\alpha-1}{2}}\right)^{2/(\alpha-1)}} \quad \Longleftrightarrow \quad B\geq \theta^{-2/(\alpha-1)},
\]
The proof of $A\leq \theta^{-2/(\alpha-1)}$ is similar.

\textbf{Step 3:} In this step we show that $A=B$,
which leads to the conclusion of the theorem.
Otherwise, we can assume $B> \theta^{-2/(\alpha-1)}$ without loss of generality.
Let $A_1$ and $B_1$ be two fixed constants such that $\theta^{-2/(\alpha-1)}< A_1<B_1 <B$.
Since $w$ is continuous and that $\liminf_{x\to\infty} w(x)<A_1<B_1< \limsup_{x\to\infty} w(x)$.  The following sequences are well-defined:
\begin{align}
	& a_1:= \inf\{x>0: w(x)= A_1\},\quad   b_1:= \inf\{x> a_1: w(x)= B_1 \}, \nonumber\\
	& a_k:= \inf\{x> b_{k-1}: w(x)= A_1\}, \quad b_k:=\inf\{x> a_k: w(x)=B_1\},\nonumber\\
	& a_k^*:= \sup\{x\in [a_k, b_k): w(x)=A_1 \}.
\end{align}
Note that $a_k\uparrow\infty$ and $b_k\uparrow\infty$.
Taking $x= a_k^*$ in \eqref{step_20},
by \eqref{uniform-limit} and noticing that $\phi(y) =(\theta y+1)^{-\frac{2}{\alpha-1}}$,
we get that for any $K>0$ and any $\varepsilon>0$ with $(1+\varepsilon)A_1< B_1$, there exists $N$ such that
\begin{align}\label{step_22}
	\sup_{y\in [0,K]}\left| \frac{w\left(a_k^*\left(1+ y A_1^{-\frac{\alpha-1}{2}}\right)\right)}{A_1} \cdot
	\frac{(\theta y+1)^{2/(\alpha-1)} }{\left( 1+ yA_1^{-\frac{\alpha-1}{2}}\right)^{2/(\alpha-1)}}-1\right|<\varepsilon,\quad k>N.
\end{align}
Since $A_1> \theta^{-2/(\alpha-1)}\Longleftrightarrow A_1^{-(\alpha-1)/2} < \theta$, by \eqref{step_22}, we see that when $k>N$,
\begin{align}
	\sup_{y\in [0,K]} w\left(a_k^*\left(1+ y A_1^{-\frac{\alpha-1}{2}}\right)\right) < (1+\varepsilon)A_1<B_1,
\end{align}
which implies that  for any $k>N$,
\begin{align}\label{step_24}
 \left\{a_k^*\left(1+y A_1^{-\frac{\alpha-1}{2}}\right): y\in [0, K] \right\}\subset [a_k^*, b_k)
\end{align}
by the definition of $b_k$.  On the other hand, for any $\delta>0$, \eqref{step_22} implies that  uniformly for all $y\in [\delta, K]$,
\begin{align}
	\lim_{k\to\infty} \frac{w\left(a_k^*\left(1+ y A_1^{-\frac{\alpha-1}{2}}\right)\right)}{A_1} <1,
\end{align}
which implies that there exists $N_1$ such that for all $k>N_1$,
\begin{align}\label{step_23}
	\sup_{y\in [\delta,K]} w\left(a_k^*\left(1+ y A_1^{-\frac{\alpha-1}{2}}\right)\right) <A_1.
\end{align}
Therefore, by the continuity of $w$ and the definitions of $a_k, b_k, a_k^*$,  for any $k>N_1$, there
exists $m_k>k$ such that
\begin{align}\label{step_25}
	\left\{a_k^*\left(1+ y A_1^{-\frac{\alpha-1}{2}}\right): y\in [\delta, K] \right\} \subset [a_{m_k}, a_{m_k}^*].
\end{align}
Moreover, for $y=K$,
\begin{align}
	a_k^*\left(1+ K A_1^{-\frac{\alpha-1}{2}}\right)\geq a_{m_k}\geq a_{k+1}> b_k,
\end{align}
which contradicts \eqref{step_24}. This completes the proof of the theorem.
\hfill$\Box$

%\bigskip
%\noindent
%{\bf Acknowledgements:}
%We thank the referees for very helpful comments on the first version of this paper.
%\bigskip
\noindent

\begin{singlespace}

\end{singlespace}

\vskip 0.2truein
\vskip 0.2truein

\noindent{\bf Haojie Hou:}  School of Mathematical Sciences, Peking
University,   Beijing, 100871, P.R. China. Email: {\texttt
houhaojie@pku.edu.cn}

\smallskip

\noindent{\bf Yiyang Jiang:}  School of Mathematical Sciences, Peking
University,   Beijing, 100871, P.R. China. Email: {\texttt
	jyy.0916@stu.pku.edu.cn}

\smallskip

\noindent{\bf Yan-Xia Ren:} LMAM School of Mathematical Sciences \& Center for
Statistical Science, Peking
University,  Beijing, 100871, P.R. China. Email: {\texttt
yxren@math.pku.edu.cn}

\smallskip
\noindent {\bf Renming Song:} Department of Mathematics,
University of Illinois,
Urbana, IL 61801, U.S.A.
Email: {\texttt rsong@illinois.edu}

\end{document}